\documentclass[12pt]{amsart}

\usepackage{amsmath}
\usepackage{amssymb}
\usepackage{url}
\usepackage{amscd}
\usepackage{mathrsfs}
\usepackage[dvips]{graphicx}
\usepackage{indentfirst}
\usepackage{setspace}

\theoremstyle{plain}

\newtheorem{lemma}{Lemma}

\newcommand{\gf}{{\mathrm{GF}}}

\newcommand{\cc}{{\mathbf c}}

\newcommand{\tr}{{\mathrm{Tr}}}

\newcommand{\F}{{\mathcal F}}

\newcommand{\FF}{{\mathbb{F}}}

\newcommand{\C}{{\mathcal C}}

\renewcommand{\vec}[1]{\underline{#1}}

\def\({\left(}
\def\){\right)}

\begin{document}

%\date{Received \today}
\title{The weight distributions of a class of cyclic codes}

%\author[A. Iosevich]{{\sc Alex Iosevich}}
%\author[I. E. Shparlinski]{\sc Igor E.~Shparlinski}
\author[M. Xiong]{\sc Maosheng Xiong}

\address{Maosheng Xiong: Department of Mathematics,
Hong Kong University of Science and Technology,
Clear Water Bay, Kowloon, Hong Kong
%% \tt iosevich@math.missouri.edu
}
\email{mamsxiong@ust.hk}

\thanks{The author was supported by the Research Grants Council of Hong Kong under Project Nos. RGC606211 and DAG11SC02.}

%\thanks{}

\begin{abstract}
Recently, the weight distributions of the duals of the cyclic codes with two zeros have been obtained for several cases in \cite{DL1,DL2,WT}. In this paper we provide a slightly different approach toward the general problem and use it to solve one more special case. We make extensive use of standard tools in number theory such as characters of finite fields, the Gauss sums and the Jacobi sums to transform the problem of finding the weight distribution into a problem of evaluating certain character sums over finite fields, which on the special case is related with counting the number of points on some elliptic curves over finite fields. Other cases are also possible by this method.
\end{abstract}

\maketitle

\thispagestyle{empty}

\maketitle

\thispagestyle{empty}
%*********************************************%
\section{Introduction}
%*********************************************%

% Let $r=q^m$ for a fixed positive integer $m$.

Denote by $\gf(q)$ the finite field of cardinality $q$, where $q=p^s$, $s$ is a positive integer and $p$ is a prime number. An $[n,k,d]$-linear code $\C$ is a $k$-dimensional subspace of $\gf(q)^n$ with minimum distance $d$. If in addition $\C$ satisfies the condition that $(c_{n-1},c_0,c_1,\ldots,c_{n-2}) \in \C$ whenever $(c_0,c_1,\ldots,c_{n-2},c_{n-1}) \in \C$, then $\C$ is called a cyclic code. Let $A_i$ denote the number of codewords with Hamming weight $i$ in $\C$. The weight enumerator of $\C$ is defined by
\[1+A_1x+A_2x^2+\cdots+A_nx^n.\]
The sequence $(1,A_1,\ldots,A_n)$ is called the weight distribution of $\C$. In coding theory it is often desirable to know the weight enumerator (or equivalently the weight distribution) of a code because they contain a lot of important information about the code, for example, they can be used to estimate the error correcting capability and the error probability of error detection and correction with respect to some algorithms. This is quite useful in practice. Many important families of cyclic codes have been studied extensively in the literature, so are their various properties, however the weight distributions are in general difficult to obtain and they are known only for a few special families. 

Given a positive integer $m$, let $r=q^m$, and $\alpha$ be a generator of $\gf(r)^*$. Let $h$ be a positive factor of $q-1$, and $e$ be a factor of $h$. Define
\begin{eqnarray} \label{1:para} g=\alpha^{(q-1)/h}, \quad n=\frac{h(r-1)}{q-1},\quad \beta=\alpha^{(r-1)/e},\quad N=\gcd\left(m,\frac{e(q-1)}{h}\right).\end{eqnarray}
The order of $g$ is $n$ and $(g \beta)^n=1$. It is also known that the minimal polynomials of $g^{-1}$ and $(\beta g)^{-1}$ are distinct over $\gf(q)$, hence their product is a divisor of $x^n-1$ (see \cite{DL1}).
Define the cyclic code
\[\C_{(q,m,h,e)}=\left\{\cc_{(a,b)}: a,b \in \gf(r)\right\},\]
where the codeword $\cc_{(a,b)}$ is give by
\[\cc_{(a,b)}:=\left(\tr\left(ag^i+b(\beta g)^i\right)\right)_{i=0}^{n-1}.\]
Here for simplicity $\tr$ is the trace function from $\gf(r)$ to $\gf(q)$.

When $h=q-1$, the code $\C_{(q,m,h,e)}$ is the dual of the primitive cyclic linear code with two zeros, which have been well studied (see for example \cite{BM, CCD, CCZ,C,Mc,MR,S,YCD}); In general the dimension of the code $\C_{(q,m,h,e)}$ is a factor of $2m$ and its weight distribution could be very complex. The known results on the weight distribution of $\C_{(q,m,h,e)}$ are listed as follows:
\begin{itemize}
\item[1)] $e>1$ and $N=1$ (\cite{DL1});

\item[2)] $e=2$ and $N=2$ (\cite{DL1});

\item[3)] $e=2$ and $N=3$ (\cite{DL2});

\item[4)] $e=2$ and $p^j+1 \equiv 0 \pmod{N}$, where $j$ is a positive integer (\cite{DL2});

\item[5)] $e=3$ and $N=2$ (\cite{WT}).

\end{itemize}
The purpose of this paper is to compute the weight distribution for one more case, that is for $e=4, N=2$. As was indicated in \cite{WT} and also in several other papers (\cite{LF1,LTW}), the problem of computing the weight distribution of a code often boils down to evaluating certain character sums or counting the number of points on a curve over a finite field. The strategy is similar, however, our treatment is different from \cite{WT} and \cite{DL1,DL2}. Our method, using orthogonality properties of characters of finite fields, allows us to translate the problem directly into the problem of evaluating certain character sums. This is still difficult in general, but on the special case that $e=4, N=2$, such character sums can be evaluated by counting the number of points on some elliptic curves over a finite field, which are fortunately well-known for a long time. The case $e=3, N=2$ which was resolved in \cite{WT} can also be handled easily by this method. %For other more general cases of $e$ and $N$ we can not evaluate the character sums, however, we can obtain a fairly good estimate on the $A_i$'s in the weight distribution of $\C_{(q,m,h,e)}$ by appealing to the Riemann hypothesis over finite fields.

We remark that the case $e=3,N=3$ can also be computed explicitly, however, the results are quite complicated for the following reasons: for $N=3$, the three Gaussian periods may have two or three distinct values depending on the parameters, and the choice of $g$ as a cubic power or not and the prime $p$ such that $p=2$, $p \equiv 1 \pmod{3}$ or $p \equiv 2 \pmod{3}$ all have a subtle influence on the weight distribution of the code. There are simply too many cases to consider and a lot of computation is involved. For the sake of clarity, it seems more appropriate to write a separate paper for the case $e=3,N=3$. Some other special cases may also be possible.

The paper is organized as follows: in Section 2 we use standard theory of characters of finite fields to set up our strategy, then in Section 3 we use the method to compute explicitly the weight distribution for the case $e=4,N=2$ (see Table 1--4 in Section 3). Finally in Section 4 we provide four examples by using Magma. We find that the papers \cite{DL1,DL2,WT} quite inspiring and very well-written, which we use as general references and starting points of this paper. Interested readers may refer to them for some preliminary background and other information related with this subject.

%*********************************************%
\section{A general strategy}
%*********************************************%

\subsection{A general strategy}

We set up our strategy for the general situation. It will be used throughout the paper. The parameters $p,q,r,\alpha,\ldots$ etc are from Section 1.

Denote by $C^{(N,r)}$ the subgroup of $\gf(r)^*$ generated by $\alpha^N$. The subset $C^{(N,r)} $ consists of non-zero elements which are perfect $N$-th powers in $\gf(r)$. This is also called the cyclotomic class of order $N$ in $\gf(r)$ with respect to $1$. Since $N|m,N|(q-1)$, the integer $(r-1)/(q-1)=q^{m-1}+q^{m-2}+\cdots+q+1$ is divisible by $N$, hence $\beta \in C^{(N,r)}$. It is also easy to see that $\gf(q)^* \subset C^{(N,r)}$.

For any $u \in \gf(r)$, define
\begin{eqnarray} \label{2:eta} \eta_u^{(N,r)}=\sum_{z \in C^{(N,r)}}\psi(zu),\end{eqnarray}
where $\psi$ is the canonical additive character of $\gf(r)$, which is given by $\psi(x)=\exp\left(\frac{2 \pi i}{p} \tr_p(x)\right)$, here $\tr_p$ is the trace function from $\gf(r)$ to $\gf(p)$. Obviously $\eta_0^{(N,r)}=\frac{r-1}{N}$. If $u \ne 0$, the term $\eta_u^{(N,r)}$ is called a ``Gaussian period'', a well-known important object which has been studied since Gauss. Note that our notation $C^{(N,r)}$ and $\eta_u^{(N,r)}$ differ slightly from \cite{DL1,DL2,WT}, however, they serve our purpose well. Also note that the Gaussian periods $\eta_u^{(N,r)}$, $u \ne 0$ depend only on the particular coset of $\gf(r)^*$ with respect to $C^{(N,r)}$ that $u$ belongs to, so there are $N$ such Gaussian periods.

The starting point of our computation is that, by \cite[Lemma 5]{DL2} (see also \cite{DL1,WT}), for any $(a,b) \in \gf(r)^2$, the Hamming weight of the codeword $\cc_{(a,b)}$ is equal to $n-Z(a,b)$, where
\begin{eqnarray} \label{2:z}
Z(a,b)=\frac{h(r-1)}{q(q-1)}+\frac{hN}{eq}\sum_{i=1}^{e} \eta_{(a+\beta^i b)g^i}^{(N,r)}\,.
\end{eqnarray}
Let us define for simplicity the ``modified weight'' of $\cc_{(a,b)}$ as
\begin{eqnarray} \label{2:lambda} \lambda(a,b)=\frac{hN}{eq}\sum_{i=1}^{e} \eta_{(a+\beta^i b)g^i}^{(N,r)}\,. \end{eqnarray}
It suffices to study $\lambda(a,b)$ only. From (\ref{2:lambda}) we see that the weight $\lambda(a,b)$ is always a simple linear combination of $\eta_0^{(N,r)}=(r-1)/N$ and the $N$ Gaussian periods $\eta_u^{(N,r)}$, $u \ne 0$. Our strategy is to at first find all possible values of $\lambda(a,b)$, and then for each such value, we count the number of $(a,b)$'s such that $\lambda(a,b)$ attains this value. There are two cases that we need to consider separately, depending on whether or not the term $\eta_0^{(N,r)}$ appears in the expression of $\lambda(a,b)$.

\subsection{Case 1.} Suppose $\prod_{i=1}^e\left(a+\beta^i b\right) \ne 0$. For any $c_1,\ldots,c_e \in \gf(r)^*$, we write $\vec{c}=(c_1,\ldots,c_e)$ and define
\[\F(\vec{c})=\left\{(a,b) \in \gf(r)^2: \begin{array}{ll}
\left(a+\beta^i b\right)g^i c_i \in C^{(N,r)} \,\, \forall i \end{array}
\right\}.\]
Then for any $(a,b) \in \F(\vec{c})$ we obtain
\begin{eqnarray} \label{2:z1}
\lambda(a,b)=\frac{hN}{eq}\sum_{i=1}^{e} \eta_{c_i^{-1}}^{(N,r)}.
\end{eqnarray}

Now we need to find the cardinality $\#\F(\vec{c})$ for each $\vec{c}$, which is denoted by $f(\vec{c})$. Applying the orthogonality property (\cite{IR}) 
\begin{eqnarray} \label{2:or1} \frac{1}{N}\sum_{\chi^N=\epsilon}\chi\left(x\right)=\left\{
\begin{array}{ll}
1:& \mbox{ if } x \in C^{(N,r)},\\
0:& \mbox{ if } x \not \in C^{(N,r)},
\end{array}
\right.\end{eqnarray}
where the sum is over all multiplicative characters $\chi$ of $\gf(r)^*$ such that $\chi^N=\epsilon$, $\epsilon$ is the principal character (we use the convention $\chi(0):=0$ to extend the definition of $\chi$ to $\gf(r)$), we can write $f(\vec{c})$ as
\begin{eqnarray} \label{2:ffc} f(\vec{c})=\sum_{a,b}\prod_{i=1}^e \left\{\frac{1}{N}\sum_{\chi_i^N=\epsilon}\chi_i\Bigl(\left(a+\beta^i b\right)g^ic_i\Bigr)\right\},\end{eqnarray}
where the outer sum is over all $a,b \in \gf(r)$, and the inner sum for each $i$ is over all multiplicative characters $\chi_i$ of $\gf(r)^*$ such that $\chi_i^N=\epsilon$.

Expanding the product on the right side, interchanging the order of summation, separating the sum for $a=0$ and $a \ne 0$, and applying the identity (\cite{IR})
\begin{eqnarray} \label{2:or2} \frac{1}{r-1}\sum_{x \in \gf(r)}\chi(x)=\left\{
\begin{array}{ll}
1:& \mbox{ if } \chi=\epsilon,\\
0:& \mbox{ if } \chi \ne \epsilon,
\end{array}
\right.\end{eqnarray}
we can obtain
\begin{eqnarray*}
f(\vec{c})=f_1(\vec{c})+f_2(\vec{c}),
\end{eqnarray*}
where
\[f_1(\vec{c})=\frac{r-1}{N^e}\sum_{\substack{\chi_i^N=\epsilon\\
\chi_1 \cdots \chi_e=\epsilon}} \prod_{i=1}^e \chi_i\left(g^i \beta^i c_i\right),\]
\[f_2(\vec{c})=\frac{r-1}{N^e}\sum_{\substack{\chi_i^N=\epsilon\\
\chi_1 \cdots \chi_e=\epsilon}} \prod_{i=1}^e \chi_i\left(g^i c_i\right) \sum_{b}\prod_{i=1}^e \chi_i\left(1+ \beta^i b\right).\]
For $f_1(\vec{c})$ and $f_2(\vec{c})$, in writing $\chi_e=\chi_1^{-1} \cdots \chi_{e-1}^{-1}$, then the sum is over all characters $\chi_i$, $1 \le i \le e-1$ such that $\chi_i^N=\epsilon$. Noticing that $\beta^e=1$, $\beta \in C^{(N,r)}$ and $g^e \in C^{(N,r)}$, we can simplify $f_1(\vec{c})$ and $f_2(\vec{c})$ further as
\[f_1(\vec{c})=\frac{r-1}{N^e}\sum_{\substack{\chi_i^N=\epsilon\\ \chi_1, \ldots, \chi_{e-1}}} \prod_{i=1}^{e-1} \chi_i\left(g^{i} \beta^{i} c_ic_e^{-1}\right),\]
\[f_2(\vec{c})=\frac{r-1}{N^e}\sum_{\substack{\chi_i^N=\epsilon\\ \chi_1, \ldots, \chi_{e-1}}} \prod_{i=1}^{e-1} \chi_i\left(g^{i} c_ic_e^{-1}\right) \sum_{b \ne -1}\prod_{i=1}^{e-1} \chi_i\left(\frac{1+ \beta^i b}{1+b}\right).\]
As for $f_2(\vec{c})$, make a change of the variable $1+b \to b'$ and then $b' \to b'^{-1}$, and then make up the term for $b'=0$, we find that
\[f_2(\vec{c})=\frac{r-1}{N^e}\sum_{\substack{\chi_i^N=\epsilon\\ \chi_1, \ldots, \chi_{e-1}}} \prod_{i=1}^{e-1} \chi_i\left(g^{i} c_ic_e^{-1}\right) \left( \sum_{b}\prod_{i=1}^{e-1} \chi_i\left(\beta^i+ (1-\beta^i) b\right) -\prod_{i=1}^{e-1} \chi_i\left(\beta^i\right)\right).\]
The second term in $f_2(\vec{c})$ is canceled out with $f_1(\vec{c})$. Hence combining $f_1(\vec{c})$ and $f_2(\vec{c})$ together we obtain
\begin{eqnarray} \label{2:fc}
f(\vec{c})=\frac{r-1}{N^e}\sum_{\substack{\chi_i^N=\epsilon\\
\chi_1, \ldots, \chi_{e-1}}} f_{\chi_1,\ldots,\chi_{e-1}}(\vec{c}),
\end{eqnarray}
where
\[f_{\chi_1,\ldots,\chi_{e-1}}(\vec{c})=\prod_{i=1}^{e-1} \chi_i\left(g^{i} (1-\beta^i) c_ic_e^{-1}\right) \sum_{b}\prod_{i=1}^{e-1} \chi_i\left(b+\gamma_i\right),\]
and
\begin{eqnarray} \label{2:gamma}
\gamma_i=\frac{\beta^i}{1-\beta^i}, \quad i=1,2, \ldots, e-1.
\end{eqnarray}

In general it might be difficult to evaluate $f(\vec{c})$ exactly, however, we can get a fairly good estimate. If $\chi_i \ne \epsilon$ for some $i$, $1 \le i \le e-1$, then appealing to Weil's bound on character sums over finite fields (\cite[Theorem 11.23]{IK}) we have
\[|f_{\chi_1,\ldots,\chi_{e-1}}(\vec{c})| \le (e-2)\sqrt{r}. \]
On the other hand, if $\chi_i=\epsilon$ for all $i$, then
\[f_{\epsilon,\ldots,\epsilon}(\vec{c})=\sum_{\substack{b\\
b+\gamma_i \ne 0\, \forall i}} 1=r-e+1. \]
So we obtain
\[\left|f(\vec{c})-\frac{(r-1)(r-e+1)}{N^e}\right|\le \frac{(e-2)(r-1)\sqrt{r}}{N}. \]
This shows that when $r$ is large compared with $N^e$, then for any $\vec{c}$, the simple linear combination of $\eta_{c_i^{-1}}^{(N,r)}$'s could appear in the weight $\lambda(a,b)$, and the frequency of such $(a,b)$'s for this to occur is of roughly the same amount.

\subsection{Case 2.} Now suppose that $a+\beta^t b=0$ for some $t$, $1 \le t \le e$ and $(a,b) \ne (0,0)$. We have $a=-\beta^tb \ne 0$ and $a+\beta^i b=b(\beta^i-\beta^t)$. Now from (\ref{2:lambda}) we obtain
\begin{eqnarray} \label{2:z2} \lambda\left(-\beta^t b,b\right)=\frac{hN}{eq} \left\{\frac{r-1}{N}+\sum_{\substack{i=1\\ i \ne t}}^e \eta_{bg^i(\beta^i-\beta^t)}^{(N,r)}\right\}, \quad 1 \le t \le e. \end{eqnarray}

%*********************************************%
\section{The case $e=4,N=2$}
%*********************************************%

When $e=4,N=2$, the parameters are
\[\beta=\alpha^{(r-1)/4}, \quad g=\alpha^{(q-1)/h}, \quad 2=\gcd\left(m, \frac{4(q-1)}{h}\right),\quad 4|h|(q-1). \]
Hence $q$ is odd and $\beta^4=1, \beta^2=-1$. We also know that $\beta$ and any $a \in \gf(q)^*$ are all squares in $\gf(r)$. The Gaussian periods are (\cite{MY})
\begin{eqnarray}
\eta_1^{(2,r)}=\left\{\begin{array}{ll}
\frac{-1-(-1)^{sm}\sqrt{r}}{2} & \mbox{ if } p \equiv 1 \pmod{4};\\
\frac{-1-(-i)^{sm}\sqrt{r}}{2} & \mbox{ if } p \equiv 3 \pmod{4};
\end{array}, \quad i=\sqrt{-1},
\right.
\end{eqnarray}
and $\eta_{\alpha}^{(2,r)}=-1-\eta_1^{(2,r)}$, $\eta_0^{(2,r)}=(r-1)/2$. %For any $(a,b) \in \gf(r)^2$, the value $\lambda(a,b)$ is always a linear combination of $\eta_{0}^{(2,r)},\eta_{1}^{(2,r)}$ and $\eta_{\alpha}^{(2,r)}$.

\subsection{Evaluation of $f(\vec{c})$}
Denote by $\chi$ the non-trivial quadratic character of $\gf(r)^*$. The $f(\vec{c})$ given in (\ref{2:fc}) can be written explicitly as
\begin{eqnarray*}
f(\vec{c})&=&\frac{r-1}{2^4} \biggl\{f_{\epsilon,\epsilon,\epsilon}(\vec{c})+f_{\epsilon,\epsilon,\chi}(\vec{c})+
f_{\epsilon,\chi,\epsilon}(\vec{c})+f_{\chi,\epsilon,\epsilon}(\vec{c})+\biggr.\\
&& \biggl.f_{\epsilon,\chi,\chi}(\vec{c})+f_{\chi,\epsilon,\chi}(\vec{c})+
f_{\chi,\chi,\epsilon}(\vec{c})+f_{\chi,\chi,\chi}(\vec{c})
\biggr\}.
\end{eqnarray*}
We will compute each term individually. From (\ref{2:gamma}) it is easy to see that
\[\gamma_1-\gamma_3=\beta, \quad \gamma_2-\gamma_3=\gamma_1-\gamma_2=\frac{\beta}{2}. \]
Noticing that $\gamma_1,\gamma_2,\gamma_3$ are all distinct, we obtain
\[f_{\epsilon,\epsilon,\epsilon}(\vec{c})= \sum_{\substack{b \\
b+\gamma_i \ne 0 \, \forall i}}1=r-3. \]
As for $f_{\epsilon,\epsilon,\chi}(\vec{c})$ we have
\[f_{\epsilon,\epsilon,\chi}(\vec{c})=\chi\bigl(g^3(1-\beta^3)c_3c_4^{-1}\bigr)
\sum_{\substack{b\\
b+\gamma_i \ne 0 \, \forall i}} \chi(b+\gamma_3). \]
Since $\chi$ is a quadratic character and
\begin{eqnarray} \label{3:beta} \chi(a \beta)=1, \quad \forall a \in \gf(q)^*, \end{eqnarray}
we have
\[\sum_{\substack{b\\
b+\gamma_i \ne 0 \, \forall i}} \chi(b+\gamma_3)= \sum_{\substack{b}} \chi(b+\gamma_3)-\chi\left(\gamma_3-\gamma_1\right)-\chi\left(\gamma_3-\gamma_2\right)=-2.\]
Using (\ref{3:beta}) again and that $\beta^2=-1$ we obtain
\[f_{\epsilon,\epsilon,\chi}(\vec{c})=-2 \chi\bigl(g(1+\beta)c_3c_4\bigr).\]
Similarly we obtain
\[f_{\epsilon,\chi,\epsilon}(\vec{c})=-2 \chi\left(c_2c_4\right),\quad
f_{\chi,\epsilon,\epsilon}(\vec{c})=-2 \chi\bigl(g(1+\beta)c_1c_4\bigr).\]
Now we compute $f_{\epsilon,\chi,\chi}(\vec{c})$. We obtain
\[f_{\epsilon,\chi,\chi}(\vec{c})=\chi\bigl(g^2(1-\beta^2)c_2c_4^{-1}\bigr)
\chi\bigl(g^3(1-\beta^3)c_3c_4^{-1}\bigr)\sum_{\substack{b\\ b+\gamma_1 \ne 0}} \chi\bigl((b+\gamma_2)(b+\gamma_3)\bigr). \]
Making up the term for $b+\gamma_1=0$, this can be simplified as
\[f_{\epsilon,\chi,\chi}(\vec{c})=\chi\bigl(g(1+\beta)c_2c_3\bigr)
\left(-1+\sum_{\substack{b}} \chi\bigl((b+\gamma_2)(b+\gamma_3)\bigr)\right). \]
Applying the following lemma we find easily that
\[f_{\epsilon,\chi,\chi}(\vec{c})=-2\chi\bigl(g(1+\beta)c_2c_3\bigr). \]

\begin{lemma} \label{3:lem1} Let $\chi$ be a non-trivial quadratic character of $\gf(r)$, and $a,b \in \gf(r)$ are distinct. Then
\[ \sum_{\substack{x \in \gf(r)}} \chi\bigl((x+a)(x+b)\bigr)=-1. \]
\end{lemma}
\noindent {\bf Proof.} Denote by $A$ the number of $\gf(r)$-rational points $(x,y,z)$ on the curve
\begin{eqnarray} \label{3:curve} (x+ay)(x+by)=z^2. \end{eqnarray}
It is known that $A$ can be written as the character sum
\[A=r^2+\sum_{x,y \in \gf(r)}\chi\left((x+ay)(x+by)\right).\]
In the sum over $x,y \in \gf(r)$, separating the cases that $y=0$ and $y \ne 0$ and then making a change of variable $x \to xy$, we find that
\begin{eqnarray} \label{3:ch} A=r^2+(r-1)+(r-1) \sum_{x \in \gf(r)}\chi\left((x+a)(x+b)\right).\end{eqnarray}
On the other hand, we can compute $A$ directly by solving the equation (\ref{3:curve}): make a change of variables
\[x+ay=\lambda, \quad x+by=\mu.\]
Since
\[y=(\lambda-\mu)/(a-b), \quad x=\lambda-ay,\]
we see that $(x,y) \leftrightarrow (\lambda,\mu)$ is a one-to-one correspondence. So $A$ is also the number of $\gf(r)$-rational points $(\lambda,\mu,z)$ on the curve
\[\lambda \mu=z^2.\]
It is easy to count that $A=r^2$. Combining this with (\ref{3:ch}) completes the proof of Lemma \ref{3:lem1}. \quad $\square$

Similarly we obtain
\[f_{\chi,\epsilon,\chi}(\vec{c})=-2\chi\left(c_1c_3\right), \quad f_{\chi, \chi,\epsilon}(\vec{c})=-2\chi\bigl(g(1+\beta)c_1c_2\bigr). \]
Finally, we need to evaluate $f_{\chi,\chi,\chi}(\vec{c})$. We have
\[f_{\chi,\chi,\chi}(\vec{c})=\chi\left(c_1c_2c_3c_4\right)\sum_{b} \chi\bigl((b+\gamma_1)(b+\gamma_2)(b+\gamma_3)\bigr).\]
Denote by $A$ the number of $\gf(r)$-rational points $(x,y)$ on the elliptic curve
\begin{eqnarray} \label{3:e2} y^2=(x+\gamma_1)(x+\gamma_2)(x+\gamma_3).\end{eqnarray}
Clearly
\[\sum_{b} \chi\bigl((b+\gamma_1)(b+\gamma_2)(b+\gamma_3)\bigr)=A-r. \]
To count $A$, we make a change of variable $x'=x+\gamma_2$, and then $x''=x/2^2,\,\,y''=y/2^3$, the curve (\ref{3:e2}) is transformed into
\begin{eqnarray} \label{3:e} E: y^2=x^3+{4x}. \end{eqnarray}
The elliptic curve $E$ is well-known, its theory and properties have been extensively studied. For example, \cite[Theorem, p. 59]{KO} computed explicitly the Zeta function of the curve $y^2=x^3-n^2x$ over any finite field $\FF_p$ (see also \cite[Exercise 23, p. 64]{KO}), and \cite[Theorem 4, p. 305]{IR} found explicitly the number of points on the curve $y^2=X^3+Dx$ over any finite field $\FF_p$, where $p$ is a prime number, $n,D$ are any integers. The proofs use standard techniques involving Gauss sums and Jacobi sums. The number of $\gf(r)$-points on the curve $E$ can also be obtained in a very similar way, with slightest modifications in the argument. Interested readers may refer to the two references for details. We record the result as follows.

\begin{lemma} \label{3:le2} Let $E$ be the elliptic curve given in (\ref{3:e}) over the finite field $\gf(p)$, where $p$ is an odd prime. For any positive integer $n \ge 1$, denote by $N_n$ the number of $\gf(p^n)$-rational points on $E$ (including the point at infinity). Then
\begin{eqnarray*} N_n=1+p^n-\pi^n-\bar{\pi}^n,\end{eqnarray*}
where $\pi$ can be computed as follows.
\begin{itemize}
\item[1).] If $p \equiv 1 \pmod{4}$, then $\pi$ is any Gaussian integer of norm $p$ such that $\pi \equiv 1 \pmod{2+2i}$.

\item[2).] If $p \equiv 3 \pmod{4}$, then $\pi=i \sqrt{p}$.
\end{itemize}
%Here $i=\sqrt{-1}$ is the complex number.
\end{lemma}

Using Lemma \ref{3:le2}, since $r=q^m=p^{ms}$, we have
\[A=r-\pi^{ms}-\bar{\pi}^{ms},\]
and therefore we can obtain
\[f_{\chi,\chi,\chi}(\vec{c})=-\chi(c_1c_2c_3c_4)\left(\pi^{ms}+\bar{\pi}^{ms}\right). \]
In summary we obtain
\begin{lemma} \label{3:le3}
For any $\vec{c}=(c_1,\ldots,c_e)$ where $c_1,\ldots,c_e \in \gf(r)^*$, we have
\begin{eqnarray*}
f(\vec{c})&=&\frac{r-1}{2^4} \biggl( r-3-2\chi\left(c_2c_4\right)-2\chi\left(c_1c_3\right)-\chi(c_1c_2c_3c_4)\left(\pi^{ms}+\bar{\pi}^{ms}\right)\biggr.\\
&&-2\chi\left(g(\beta+1)\right)\Bigl\{\chi(c_3c_4)+\chi(c_1c_4)+\chi\left(c_2c_3\right)+\chi\left(c_1c_2\right)\Bigr\}\biggr).
\end{eqnarray*}
\end{lemma}

\subsection{The remaining case}
Next we need to evaluate $\lambda\left(-\beta^tb,b\right)$ in (\ref{2:z2}). We adopt a notation: $\lambda \equiv \mu \pmod{\square}$ means that $\lambda \mu \in \gf(r)^*$ is a square; $\lambda=\square$ means that $\lambda \in \gf(r)^*$ is a square.

For $t=1$, it is easy to see that
\[\beta^2-\beta \equiv \beta^4-\beta \equiv \beta+1 \pmod{\square}, \quad \beta^3-\beta \equiv 1 \pmod{\square},\]
so we obtain
\begin{eqnarray} \label{3:r1}\lambda\left(-\beta b,b\right)=\frac{2h}{4q} \left\{\frac{r-1}{2}+2 \eta_{b(\beta+1)}^{(2,r)}+\eta_{bg}^{(2,r)}\right\}.\end{eqnarray}

Similarly we find
\begin{eqnarray} \label{3:r2} \lambda\left(-\beta^2 b,b\right)=\frac{2h}{4q} \left\{\frac{r-1}{2}+2 \eta_{bg(\beta+1)}^{(2,r)}+\eta_{b}^{(2,r)}\right\}, \quad t=2,\end{eqnarray}
\begin{eqnarray} \label{3:r3}\lambda\left(-\beta^3 b,b\right)=\frac{2h}{4q} \left\{\frac{r-1}{2}+2 \eta_{b(\beta+1)}^{(2,r)}+\eta_{bg}^{(2,r)}\right\},\quad t=3, \end{eqnarray}
\begin{eqnarray} \label{3:r4}\lambda\left(-\beta b,b\right)=\frac{2h}{4q} \left\{\frac{r-1}{2}+2 \eta_{bg(\beta+1)}^{(2,r)}+\eta_{b}^{(2,r)}\right\}, \quad t=4.\end{eqnarray}

\subsection{Conclusion}
For any $(a,b) \in \F(\vec{c})$ we have
\begin{eqnarray*}
\lambda(a,b)=\frac{2h}{4q}\sum_{i=1}^{4} \eta_{c_i^{-1}}^{(2,r)}.
\end{eqnarray*}
Since $u \ne 0$,
\begin{eqnarray} \label{3:g} \eta_u^{(2,r)}=\left\{\begin{array}{ll}
\eta_1^{(2,r)} & \mbox{ if } u=\square, \\
\eta_\alpha^{(2,r)} & \mbox{ if } u \ne \square, %(r-1)/2 & \mbox{ if } u=0.
\end{array} \right.\end{eqnarray}
Each $c_i (\ne 0)$ has only two values: either $c_i=\square$ or $c_i \ne \square$, so $\lambda(a,b)$ has at most 16 different values, and for each such value, $\#\F(\vec{c})=f(\vec{c})$ is given by Lemma \ref{3:le3}. It is clear that the result would depend on whether or not $g(\beta+1)=\square$.

Assume first that $g(\beta+1)=\square$. Then by Lemma \ref{3:le3} we have
\begin{eqnarray*}
f(\vec{c})&=&\frac{r-1}{2^4} \biggl( r-3-2\chi\left(c_2c_4\right)-2\chi\left(c_1c_3\right)-\chi(c_1c_2c_3c_4)\left(\pi^{ms}+\bar{\pi}^{ms}\right)\biggr.\\
&&-2\chi(c_3c_4)-2\chi(c_1c_4)-2\chi\left(c_2c_3\right)-2\chi\left(c_1c_2\right)\biggr).
\end{eqnarray*}
If $c_1=c_2=c_3=c_4=\square$, then the ``modified weight'' $\lambda(a,b)= \frac{2h}{4q}4 \eta_1^{(2,r)}= \frac{2h\eta_1^{(2,r)}}{q}$, and the total number of such $(a,b)$'s counted in this $\F(\vec{c})$ is given by
\[f(\vec{c})=\frac{r-1}{2^4} \left(r-15-\pi^{ms}-\bar{\pi}^{ms}\right).\]
If say $c_1=c_2=c_3=c_4 \ne \square$, then $\lambda(a,b)= \frac{2h}{4q}4 \eta_{\alpha}^{(2,r)}= \frac{2h\eta_{\alpha}^{(2,r)}}{q}$, and the total number of such $(a,b)$'s counted in this $\F(\vec{c})$ is also given by
\[f(\vec{c})=\frac{r-1}{2^4} \left(r-15-\pi^{ms}-\bar{\pi}^{ms}\right).\]
We can calculate for the other 14 cases of $\vec{c}$ in a similar way and we summarize the result in {\bf Table 1}.

\begin{table}[ht]
\caption{Part 1: The modified weight distribution for $e=2,N=4$ on the case $g(\beta+1)=\square$}
\centering
\begin{tabular}{|c| c|}
\hline
\hline
%inserts double horizontal lines
Weight $\lambda(a,b)$ & Frequency  \\
[0.5ex]
\hline
% inserts single horizontal line
$2h \eta_1^{(2,r)}/q$ & $(r-1)\left(r-15-\pi^{ms}-\bar{\pi}^{ms}\right)/16$ \\
\hline
% inserting body of the table
$2h \eta_{\alpha}^{(2,r)}/q$ & $(r-1)\left(r-15-\pi^{ms}-\bar{\pi}^{ms}\right)/16$ \\
\hline
$2h \left(3\eta_1^{(2,r)}+\eta_{\alpha}^{(2,r)}\right)/4q$ & $(r-1)\left(r-3+\pi^{ms}+\bar{\pi}^{ms}\right)/4$ \\
\hline
$2h \left(\eta_1^{(2,r)}+3\eta_{\alpha}^{(2,r)}\right)/4q$ & $(r-1)\left(r-3+\pi^{ms}+\bar{\pi}^{ms}\right)/4$ \\
\hline
$-h /q$ & $3(r-1)\left(r+1-\pi^{ms}-\bar{\pi}^{ms}\right)/8$ \\
[1ex]
% [1ex] adds vertical space
\hline
%inserts single line
\end{tabular}
\label{table:11}
\end{table}

If $g(\beta+1) \ne \square$, then
\begin{eqnarray*}
f(\vec{c})&=&\frac{r-1}{2^4} \biggl( r-3-2\chi\left(c_2c_4\right)-2\chi\left(c_1c_3\right)-\chi(c_1c_2c_3c_4)\left(\pi^{ms}+\bar{\pi}^{ms}\right)\biggr.\\
&&+2\chi(c_3c_4)+2\chi(c_1c_4)+2\chi\left(c_2c_3\right)+2\chi\left(c_1c_2\right)\biggr).
\end{eqnarray*}
The possible modified weight $\lambda(a,b)$ that could appear is the same as in Table 1, but the number of the $(a,b)$'s that could attain such weight is different. We do a similar analysis and summarize the result in {\bf Table 2}.

\begin{table}[ht]
\caption{Part 1: The modified weight distribution for $e=2,N=4$ on the case $g(\beta+1) \ne \square$}
\centering
\begin{tabular}{|c| c|}
\hline
\hline
%inserts double horizontal lines
Weight $\lambda(a,b)$ & Frequency  \\
[0.5ex]
\hline
% inserts single horizontal line
$2h \eta_1^{(2,r)}/q$ & $(r-1)\left(r+1-\pi^{ms}-\bar{\pi}^{ms}\right)/16$ \\
\hline
% inserting body of the table
$2h \eta_{\alpha}^{(2,r)}/q$ & $(r-1)\left(r+1-\pi^{ms}-\bar{\pi}^{ms}\right)/16$ \\
\hline
$2h \left(3\eta_1^{(2,r)}+\eta_{\alpha}^{(2,r)}\right)/4q$ & $(r-1)\left(r-3+\pi^{ms}+\bar{\pi}^{ms}\right)/4$ \\
\hline
$2h \left(\eta_1^{(2,r)}+3\eta_{\alpha}^{(2,r)}\right)/4q$ & $(r-1)\left(r-3+\pi^{ms}+\bar{\pi}^{ms}\right)/4$ \\
\hline
$-h /q$ & $(r-1)\left(3r-13-3\pi^{ms}-3\bar{\pi}^{ms}\right)/8$ \\
[1ex]
% [1ex] adds vertical space
\hline
%inserts single line
\end{tabular}
\label{table:12}
\end{table}

Next we need to count the frequency of the $(a,b)$'s such that they attain the modified weight $\lambda\left(-\beta^t b,b\right)$ for some $t$, $1 \le t \le 4$ which appears in (\ref{3:r1})--(\ref{3:r4}). The results also depend on whether or not $g(\beta+1)=\square$.

Assume first that $g(\beta+1)=\square$. Considering $\lambda\left(-\beta b,b\right)$ in (\ref{3:r1}), we find that
\[ \lambda\left(-\beta b,b\right)=\frac{2h}{4q}\left\{\frac{r-1}{2}+3 \eta_{bg}^{(2,r)}\right\}.\]
This is either $\frac{2h}{4q}\left\{\frac{r-1}{2}+3 \eta_{1}^{(2,r)}\right\}$ if $bg=\square$, and the number of such $b$'s is $(r-1)/2$, or $\frac{2h}{4q}\left\{\frac{r-1}{2}+3 \eta_{\alpha}^{(2,r)}\right\}$ if
$bg \ne \square$, and the number of such $b$'s is also $(r-1)/2$. It turns out that the $\lambda\left(-\beta^t b,b\right)$'s, $2 \le t \le 4$ all follow the same pattern. We summarize the result in {\bf Table 3}.

\begin{table}[ht]
\caption{Part 2: The modified weight distribution for $e=2,N=4$ on the case $g(\beta+1) = \square$}
\centering
\begin{tabular}{|c| c|}
\hline
\hline
%inserts double horizontal lines
Weight $\lambda(a,b)$ & Frequency  \\
[0.5ex]
\hline
% inserts single horizontal line
$\frac{h}{2q}\left(\frac{r-1}{2}+3 \eta_{1}^{(2,r)}\right)$ & $2(r-1)$ \\
\hline
% inserting body of the table
$\frac{h}{2q}\left(\frac{r-1}{2}+3 \eta_{\alpha}^{(2,r)}\right)$ & $2(r-1)$ \\
[1ex]
\hline
%inserts single line
\end{tabular}
\label{table:21}
\end{table}

The argument for the case that $g(\beta+1) \ne \square$ is similar. It turns out that the frequency is the same, but the weight achieved by those $(a,b)$'s is different. We record the result in {\bf Table 4}.
\begin{table}[ht]
\caption{Part 2: The modified weight distribution for $e=2,N=4$ on the case $g(\beta+1) \ne \square$}
\centering
\begin{tabular}{|c| c|}
\hline
\hline
%inserts double horizontal lines
Weight $\lambda(a,b)$ & Frequency  \\
[0.5ex]
\hline
% inserts single horizontal line
$\frac{h}{2q}\left(\frac{r-3}{2}+ \eta_{1}^{(2,r)}\right)$ & $2(r-1)$ \\
\hline
% inserting body of the table
$\frac{h}{2q}\left(\frac{r-3}{2}+ \eta_{\alpha}^{(2,r)}\right)$ & $2(r-1)$ \\
[1ex]
\hline
%inserts single line
\end{tabular}
\label{table:22}
\end{table}

Notice that the weight of $\cc_{(a,b)}$ is given by
\[w\left(\cc_{(a,b)}\right)= n- \frac{h(r-1)}{q(q-1)}-\lambda(a,b). \]
If $g(\beta+1)=\square$ then Table 1 and Table 3 (or if $g(\beta+1) \ne \square$, then Table 2 and Table 4 respectively) with the extra point corresponding to $(a,b)=(0,0)$ give the weight distribution of $\C_{q,m,h,e}$. There are always eight different weight in the code.

\section{Examples} 

\noindent {\bf Example 1.} Let $p=q=17,s=1,m=2,e=h=4,N=2$. Since $17 \equiv 1 \pmod{4}$ and $17=1^2+4^2$, we may choose $\pi=1+4i$. The Gaussian periods are $\eta_1^{(2,r)}=-9, \eta_{\alpha}^{(2,r)}=8$. Let $\alpha$ be the generator of $\gf(17^2)$ constructed from Magma, then $g(\beta+1)=\alpha^{166}$ is a square, so we may use Table 1 and Table 3 to find that the weight distribution of the cyclic code $\C_{(q,m,h,e)}$ is
\[1+576x^{48}+576x^{54}+5472x^{64}+18432x^{66}+34560x^{68}+18432x^{70}+5472x^{72}. \]
This is a $[72,4,48]$-cyclic code over $\gf(17^2)$.

\noindent {\bf Example 2.} Let $p=q=13,s=1,m=2,e=h=4,N=2$. Since $13 \equiv 1 \pmod{4}$ and $13=3^2+2^2$, we may choose $\pi=3+2i$. The Gaussian periods are $\eta_1^{(2,r)}=-7, \eta_{\alpha}^{(2,r)}=6$. Let $\alpha$ be the generator of $\gf(13^2)$ constructed from Magma, then $g(\beta+1)=\alpha^{115}$ is not a square, so we may use Table 2 and Table 4 to find that the weight distribution of the cyclic code $\C_{(q,m,h,e)}$ is
\[1+336x^{38}+336x^{40}+1680x^{48}+7392x^{50}+9744x^{52}+7392x^{54}+1680x^{56}. \]
This is a $[56,4,38]$-cyclic code over $\gf(13^2)$.

\noindent {\bf Example 3.} Let $p=3,q=3^2,s=2,m=2,e=h=4,N=2$. Since $p \equiv 3 \pmod{4}$, we have $\pi=\sqrt{3}i$. The Gaussian periods are $\eta_1^{(2,r)}=-5, \eta_{\alpha}^{(2,r)}=4$. Let $\alpha$ be the generator of $\gf(3^4)$ constructed from Magma, then $g(\beta+1)=\alpha^{72}$ is a square, so we may use Table 1 and Table 3 to find that the weight distribution of the cyclic code $\C_{(q,m,h,e)}$ is
\[1+160x^{24}+160x^{30}+240x^{32}+1920x^{34}+1920x^{36}+1920x^{38}+240x^{40}. \]
This is a $[40,4,24]$-cyclic code over $\gf(3^4)$.

\noindent {\bf Example 4.} The parameters are the same as in Example 3, except that we choose $h=8$. Then $g(\beta+1)=\alpha^{71}$ is not a square, so we may use Table 2 and Table 4 to find that the weight distribution of the cyclic code $\C_{(q,m,h,e)}$ is
\[1+160x^{52}+160x^{56}+320x^{64}+1920x^{68}+1760x^{72}+1920x^{76}+320x^{80}. \]
This is a $[80,4,52]$-cyclic code over $\gf(3^4)$.

%\[\]

%\noindent {\bf Remark.} For the case $e=5,N=2$, the method works nicely except possibly for the last term
%\[f_{\chi,\chi,\chi,\chi}(\vec{c})= \prod_{i=1}^{4} \chi\left(g^{i} (1-\beta^i) c_ic_5^{-1}\right) \sum_{b}\prod_{i=1}^{4} \chi\left(b+\frac{\beta^i}{1-\beta^i}\right).\]
%The evaluation is related with counting the number of $\gf(r)$-rational points $(x,y)$ on the curve
%\[E: y^2=\prod_{i=1}^{4} \left(x+\frac{\beta^i}{1-\beta^i}\right). \]
%Notice that this is also an elliptic curve. It is possible that exact formulas can be found, so that $f_{\chi,\chi,\chi,\chi}(\vec{c})$ can be evaluated. However, it seems quite complicated, because in order to find the weight distribution function, we also need to find exact formulas of the number of $\gf(r)$-rational points on three other elliptic curves.

%\section*{Acknowledgements}

\end{document}